\documentclass{elsart}
\usepackage{amsmath}
\usepackage{amssymb}

\newcommand{\N}{\mathbb{N}}                     
\newcommand{\R}{\mathbb{R}}                     

\newcommand{\BOne}{\mathbf{1}}
\newcommand{\bx}{\mathbf{x}}

\begin{document}

\begin{frontmatter}

\title{Approximations of Lov\'asz extensions and their induced interaction index}

\author{Jean-Luc Marichal} \ead{jean-luc.marichal[at]uni.lu}
\address{
Applied Mathematics Unit, University of Luxembourg\\
162A, avenue de la Fa\"{\i}encerie, L-1511 Luxembourg, Luxembourg }

\author{Pierre Mathonet} \ead{p.mathonet[at]ulg.ac.be}
\address{
Department of Mathematics, University of Li\`ege\\
Grande Traverse 12, B37, B-4000 Li\`ege, Belgium }

\date{June 26, 2007}

\begin{abstract}
The Lov\'asz extension of a pseudo-Boolean function $f : \{0,1\}^n \to \R$ is defined on each simplex of the standard triangulation of $[0,1]^n$
as the unique affine function $\hat f : [0,1]^n \to \R$ that interpolates $f$ at the $n+1$ vertices of the simplex. Its degree is that of the
unique multilinear polynomial that expresses $f$. In this paper we investigate the least squares approximation problem of an arbitrary Lov\'asz
extension $\hat f$ by Lov\'asz extensions of (at most) a specified degree. We derive explicit expressions of these approximations. The
corresponding approximation problem for pseudo-Boolean functions was investigated by Hammer and Holzman \cite{HamHol92} and then solved
explicitly by Grabisch, Marichal, and Roubens~\cite{GraMarRou00}, giving rise to an alternative definition of Banzhaf interaction index.
Similarly we introduce a new interaction index from approximations of $\hat f$ and we present some of its properties. It turns out that its
corresponding power index identifies with the power index introduced by Grabisch and Labreuche \cite{GraLab01}.
\end{abstract}

\begin{keyword}
pseudo-Boolean function \sep Lov\'asz extension \sep discrete Choquet integral \sep least squares approximation \sep interaction index.
\end{keyword}
\end{frontmatter}

\section{Introduction and problem setting}

Real valued set functions are extensively used in both cooperative game theory and multicriteria decision making (see for instance
\cite{GraMarRou00} and the references therein). Indeed, in its characteristic form, a {\em cooperative game}\/ on a finite set of players
$N:=\{1,\ldots,n\}$, is a set function $v:2^N\to\R$ such that $v(\varnothing)=0$ which assigns to each coalition $S$ of players a real number
$v(S)$ representing the {\em worth}\/ of $S$. In multicriteria decision making, when $N$ represents a set of criteria, a {\em capacity}\/ (also
called {\em non-additive measure}) on $N$ is a monotone set function $v:2^N\to\R$ such that $v(\varnothing)=0$ and $v(N)=1$ which assigns to
each subset $S$ of criteria its weight $v(S)$.

There is a natural identification of real valued set functions $v:2^N\to\R$ with pseudo-Boolean functions $f : \{0,1\}^n \to \R$ through the
equality $v(S)=f(\BOne_S)$, where $\BOne_S$ denotes the characteristic vector of $S$ in $\{0,1\}^n$ (whose $i$th component is $1$ if and only if
$i\in S$). Moreover, it is known \cite{HamRud68} that any pseudo-Boolean function $f : \{0,1\}^n \to \R$ has a unique expression as a
multilinear polynomial in $n$ variables
\begin{equation}\label{eq:linexprf}
f(\bx) = \sum_{S \subseteq N} a(S) \, \prod_{i \in S} x_i,
\end{equation}
where the set function $a:2^N\to\R$ is simply the M\"obius transform (see for instance \cite{Rot64}) of $v$, which can be calculated through
\begin{equation}\label{eq:mo}
a(S) = \sum_{T\subseteq S} (-1)^{|S|-|T|}\, v(T).
\end{equation}

Let $\mathfrak{S}_n$ denote the family of all permutations $\sigma$ on $N$. The {\em Lov\'asz extension}\/ $\hat f : [0,1]^n \to \R$ of any
pseudo-Boolean function $f$ is defined on each $n$-simplex
$$
S_\sigma := \{\bx \in [0,1]^n \mid x_{\sigma(1)} \leqslant \cdots \leqslant x_{\sigma(n)}\} \qquad (\sigma \in \mathfrak{S}_n),
$$
as the unique affine function which interpolates $f$ at the $n+1$ vertices of $S_\sigma$ (see Lov\'asz~\cite[\S 3]{Lov83} and Singer~\cite[\S
2]{Sin84}).

The concept of Lov\'asz extension was initially used in combinatorial optimization \cite{Fuj05,Lov83,Sin84} and then in cooperative game theory
\cite{AlgBilFerJim04,GraMarRou00}. In this paper we are mostly concerned by the interpretation of this concept in decision making, where it was
observed in \cite{Mar00f} that the Lov\'asz extension of a monotone pseudo-Boolean function is nothing but a discrete {\em Choquet integral}\/
\cite{Cho53,Den94}. Recall that discrete Choquet integrals were recently proposed in decision making as aggregation functions that generalize
the weighted arithmetic means by the taking into consideration of the interaction phenomena among input variables; see for instance
\cite{GraMurSug00,Mar02b}.

As pointed out in Grabisch et al.~\cite{GraMarRou00}, the Lov\'asz extension of any pseudo-Boolean function $f$ can be expressed as the
following min-polynomial
$$
\hat f(\bx) = \sum_{S \subseteq N} a(S) \, \min_{i \in S} x_i.
$$
Its degree is that of the multilinear polynomial (\ref{eq:linexprf}) that expresses $f$.

Hammer and Holzman~\cite{HamHol92} investigated the approximation of a pseudo-Boolean function $f$ by a multilinear polynomial of (at most) a
specified degree. More precisely, fixing $k \in \N$, with $k \leqslant n$, they defined the best $k$th approximation of $f$ as the multilinear
polynomial $f_k : \{0,1\}^n \to \R$ of degree at most $k$ which minimizes
$$
\sum_{\bx \in \{0,1\}^n} \left[f(\bx) - g(\bx)\right]^2
$$
among all multilinear polynomials $g$ of degree at most $k$. They obtained an explicit form of the solution for $k=1$ and $k=2$. For a general
$k$, the explicit form of the solution was obtained by Grabisch  et al.~\cite[\S 7]{GraMarRou00}. It is written as
$$
f_k(\bx) = \sum_{\textstyle{S \subseteq N \atop |S| \leqslant k}} a_k(S) \, \prod_{i \in S} x_i,
$$
where
\begin{equation}\label{eq:gr/ma/ro}
a_k(S) = a(S) + (-1)^{k+|S|} \sum_{\textstyle{T \supseteq S \atop |T| > k}} \frac{{|T|-|S|-1 \choose k-|S|}}{2^{|T|-|S|}} \, a(T).
\end{equation}
Hammer and Holzman~\cite{HamHol92} also noticed that $a_1(\{i\})$ identifies with the classical Banzhaf power index \cite{Ban65} related to the
$i$th player. More generally, Grabisch et al.~\cite{GraMarRou00} noticed that, for any $S\subseteq N$, the coefficient $a_{|S|}(S)$ identifies
with the Banzhaf interaction index \cite{Rou96} related to the coalition $S$.

In this paper we investigate the least squares approximation of a given Lov\'asz extension $\hat f$ by a min-polynomial of (at most) a specified
degree. More precisely, fixing $k \in \N$, with $k \leqslant n$, we give explicit formulas (see Section 2) for the best $k$th approximation of
$\hat f$, that is, the min-polynomial $\hat f_k : [0,1]^n \to \R$ of degree at most $k$ which minimizes
$$
\int_{[0,1]^n} \left[\hat f(\bx) - \hat g(\bx)\right]^2 \, d\bx
$$
among all min-polynomials $\hat g$ of degree at most $k$.

It is important to emphasize that this least squares problem has a relevant application in multicriteria decision making. Indeed, consider a
multicriteria decision making problem modelled by a Choquet integral and its underlying capacity on a set of $n$ criteria. It is sometimes
useful to approximate this Choquet integral (which is, actually, the Lov\'asz extension of the capacity) by simpler aggregation functions such
as Lov\'asz extensions of degree $k$ ($k< n$). Indeed this subclass of functions includes the Choquet integrals constructed from $k$-order
capacities introduced by Grabisch \cite{Gra97}.

By analogy with the Banzhaf interaction index obtained from the approximations of $f$, we show (see Section 3) that the approximations of $\hat
f$ give rise to a new interaction index having appealing properties (see for instance \cite{FujKojMar06}). Furthermore, it turns out that its
corresponding power index identifies with a power index introduced axiomatically by Grabisch and Labreuche \cite{GraLab01}.

In order to avoid a heavy notation, cardinalities of subsets $S,T,\ldots$ are denoted whenever possible by the corresponding lower case letters
$s,t,\ldots$, otherwise by the standard notation $|S|,|T|,\ldots$ Also, for any $S\subseteq N$ we define the function ${\rm min}_S:[0,1]^n\to\R$
as
$$
{\rm min}_S(\bx):=\min_{i \in S} x_i,
$$
with the convention that ${\rm min}_{\varnothing}\equiv 1$. Finally, we denote by $\mathcal{F}_N$ the family of set functions $v:2^N\to\R$.

\section{Approximations of Lov\'asz extensions}

We first recall the setting of our approximation problem. For any $k\in\{0,\ldots,n\}$, we denote by $V_k$ the set of all min-polynomials $\hat
f_k:[0,1]^n\to\R$ of degree at most $k$, i.e., of the form
$$
\hat f_k(\bx) = \sum_{\textstyle{S \subseteq N \atop |S| \leqslant k}} a_k(S) \, \min_{i \in S} x_i.
$$

\begin{defn}\label{de:approx}
Let $\hat f\in V_n$ and $k\in\{0,\ldots,n\}$. The {\em best $k$th approximation}\/ of $\hat f$ is the min-polynomial $\hat f_k\in V_k$ that
minimizes
$$
\int_{[0,1]^n} \left[\hat f(\bx) - \hat g(\bx)\right]^2 \, d\bx
$$
among all min-polynomials $\hat g$ of $V_k$.
\end{defn}

In this section we first discuss existence and uniqueness of the best $k$th approximation of any function $\hat f\in V_n$ (see
Proposition~\ref{exist}). Then we provide a closed-form formula for this approximation (see Theorem~\ref{thm:main}).

We begin with the following straightforward proposition:

\begin{prop}
For any $k\in\{0,\ldots,n\}$, the set $V_k$ is a linear subspace of $V_n$, which in turn is a linear subspace of the space $C_0([0,1]^n)$ of
continuous functions on $[0,1]^n$. Moreover, the set
\begin{equation}\label{eq:basisVk}
B_k := \{ {\rm min}_S \mid S \subseteq N,\; s \leqslant k \},
\end{equation}
forms a basis of $V_k$ and hence $\dim(V_k)=\sum_{s=0}^k {n \choose s}$.
\end{prop}

Consider the usual inner product $\langle \cdot \, ,\cdot \rangle$ of $C_0([0,1]^n)$ given by
$$
\langle f_1, f_2 \rangle := \int_{[0,1]^n} f_1(\bx) \, f_2(\bx) \, d\bx.
$$
This clearly induces an inner product in $V_n$, which in turn allows to define a norm $\| \hat f \| := \langle \hat f,\hat f \rangle ^{1/2}$ and
then a distance $d(\hat f_1,\hat f_2) := \|\hat f_1 - \hat f_2\|$ in $V_n$.

According to this terminology, the best $k$th approximation of $\hat f\in V_n$ is the solution $\hat f_k$ of the following least squares problem
$$
\min\{\|\hat f-\hat g\|^2 \mid \hat g\in V_k\},
$$
which amounts to a classical orthogonal projection problem.

From now on, for any $k\in\{0,\ldots,n\}$, we denote by $A_k$ the orthogonal projector from $V_n$ onto $V_k$.

\begin{prop}\label{exist}
For any $\hat f\in V_n$ and any $k\in\{0,\ldots,n\}$, the best $k$th approximation of $\hat f$ exists and is uniquely given by $\hat
f_k=A_k(\hat f)$. In particular, it is characterized as the unique element $\hat f_k$ of $V_k$ such that $\hat f - \hat f_k$ is orthogonal to
$V_k$, i.e., to each element of the basis $B_k$ defined in (\ref{eq:basisVk}).
\end{prop}

This projection problem can be simplified because of the following result, which immediately follows from the fact that the subspaces $V_k$ are
nested:

\begin{prop}\label{commute}
The operators $A_k$, $k=0,\ldots,n$, are such that
\begin{equation}\label{eq:commutation}
k \leqslant k' \quad \Rightarrow \quad A_k(A_{k'}(\hat f)) = A_k(\hat f).
\end{equation}
\end{prop}

We thus observe that $A_k(\hat f)$ can be attained from $\hat f$ by successively carrying out the projections $A_{n-1}, A_{n-2}, \ldots, A_k$.
We will therefore search for the relation that links any two consecutive projections. This relation will enable us to compute gradually
$A_k(\hat f)$ from $\hat f$.

\begin{rem}
Fix $k \in \{0,\ldots,n-1\}$ and assume that $\hat f_{k+1} = A_{k+1}(\hat f)$ is given. By (\ref{eq:commutation}), $A_k(\hat f)$ is the
orthogonal projection of $A_{k+1}(\hat f)$ onto $V_k$ and by Proposition \ref{exist}, it is given as the unique solution $\hat f_k\in V_k$ of
the following system
$$
\langle\hat f_k, {\rm min}_S\rangle = \langle\hat f_{k+1}, {\rm min}_S\rangle \qquad (S \subseteq N, \; s \leqslant k).
$$
\end{rem}

To solve this projection problem, we need the following two lemmas. The first one yields an explicit expression of the inner product $\langle
{\rm min}_S, {\rm min}_T\rangle$ of any two elements of $B_n$. The second lemma, which makes use of the first one, yields an explicit expression
of the projections onto $V_k$ of the elements of $B_{k+1}\setminus B_k$.

\begin{lem}\label{lemma:integral}
For any $S,T\subseteq N$, there holds
\begin{equation}\label{eq:int/min/min}
\int_{[0,1]^n}{\rm min}_S(\bx) \, {\rm min}_T(\bx) \, d\bx=\frac 1{|S \cup T|+2}\,\Big(\frac 1{s+1}+\frac 1{t+1}\Big).
\end{equation}
\end{lem}

\begin{pf*}{Proof.}
See Appendix A.\qed
\end{pf*}

\begin{lem}\label{lemma:enfer/1}
For any $S \subseteq N$, with $s = k+1$, we have
$$
A_k({\rm min}_{S})=\sum_{T \varsubsetneq S} (-1)^{k+t} \, \frac{{k+t+1 \choose k+1}}{{2k+2 \choose k+1}}\,{\rm min}_{T}.
$$
\end{lem}

\begin{pf*}{Proof.}
See Appendix B.\qed
\end{pf*}

\begin{prop}\label{prop:enfer}
Let $k \in \{0,\ldots,n-1\}$. Given the coefficients $a_{k+1}(S)$ of $A_{k+1}(\hat f)$, the coefficients $a_k(S)$ of $A_k(\hat f)$ are given by
\begin{equation}\label{eq:rec}
a_k(S) = a_{k+1}(S) + (-1)^{k+s} \, \frac {{k+s+1 \choose k+1}}{{2k+2 \choose k+1}} \sum_{\textstyle{T \supseteq S \atop t=k+1}} a_{k+1}(T).
\end{equation}
\end{prop}

\begin{pf*}{Proof.}
Using formula~(\ref{eq:commutation}) and then Lemma~\ref{lemma:enfer/1}, we obtain
\begin{eqnarray*}
A_k(\hat f) &=& A_k(A_{k+1}(\hat f))= \sum_{\textstyle{S \subseteq N \atop s \leqslant k+1}}a_{k+1}(S)\,A_k({\rm min}_{S})\\
&=& \sum_{\textstyle{S \subseteq N \atop s \leqslant k}}a_{k+1}(S)\,{\rm min}_{S}+\sum_{\textstyle{S \subseteq N \atop
s= k+1}} a_{k+1}(S)\,A_k({\rm min}_{S})\\
&=& \sum_{\textstyle{S \subseteq N \atop s \leqslant k}}a_{k+1}(S)\,{\rm min}_{S}+\sum_{\textstyle{S \subseteq N \atop s= k+1}}a_{k+1}(S)\sum_{T
\varsubsetneq S} (-1)^{k+t} \, \frac{{k+t+1 \choose k+1}}{{2k+2 \choose k+1}}\,{\rm min}_{T},
\end{eqnarray*}
which proves the result.\qed
\end{pf*}

We are now ready to derive the explicit form of the coefficients of $A_k(\hat f)$ in terms of the coefficients of $\hat f$. It is worth
comparing this formula with the solution (\ref{eq:gr/ma/ro}) of the Hammer-Holzman approximation problem.

\begin{thm}\label{thm:main}
Let $k \in \{0,\ldots,n\}$. The coefficients of $A_k(\hat f)$ are given from those of $\hat f$ by
\begin{equation}\label{eq:solution/system}
a_k(S) = a(S) + (-1)^{k+s} \sum_{\textstyle{T \supseteq S \atop t > k}} \frac{{k+s+1 \choose k+1}{t-s-1 \choose k-s}}{{k+t+1 \choose k+1}} \,
a(T).
\end{equation}
\end{thm}

\begin{pf*}{Proof.}
We only have to prove that the coefficients $a_k(S)$ given in (\ref{eq:solution/system}) fulfill the recurrence relation (\ref{eq:rec}).

Let $k \in \{0,\ldots,n\}$ and fix $S \subseteq N$ with $s \leqslant k$. By substituting (\ref{eq:solution/system}) into (\ref{eq:rec}), we
obtain, after removing the common term $a(S)$,
\begin{eqnarray*}
\lefteqn{(-1)^{k+s} \sum_{\textstyle{R \supseteq S \atop r > k}} \frac{{k+s+1 \choose k+1}{r-s-1 \choose k-s}}{{k+r+1
\choose k+1}} \, a(R)} \\
& = & (-1)^{k+s+1} \sum_{\textstyle{R \supseteq S \atop r > k+1}} \frac{{k+s+2 \choose k+2}{r-s-1 \choose
k-s+1}}{{k+r+2 \choose k+2}} \, a(R) \\
& & \null +(-1)^{k+s} \, \frac{{k+s+1 \choose k+1}}{{2k+2 \choose k+1}} \, \sum_{\textstyle{T \supseteq S \atop t=k+1}} \biggl[a(T) +
\sum_{\textstyle{R \supseteq T \atop r > k+1}} \frac{{2k+3 \choose k+2}}{{k+r+2 \choose k+2}} \, a(R) \biggr].
\end{eqnarray*}
Let us show that this equality holds. Dividing through by $(-1)^{k+s}$ and then removing the terms corresponding to $a(R)$, with $r = k+1$, the
equality becomes
\begin{eqnarray*}
\lefteqn{\sum_{\textstyle{R \supseteq S \atop r > k+1}} \frac{{k+s+1 \choose k+1}{r-s-1 \choose k-s}}{{k+r+1 \choose k+1}} \, a(R)}\\
&=& \null - \sum_{\textstyle{R \supseteq S \atop r > k+1}} \frac{{k+s+2 \choose k+2}{r-s-1 \choose k-s+1}}{{k+r+2 \choose k+2}} \, a(R) +
\sum_{\textstyle{R \supseteq S \atop r > k+1}} \frac{{k+s+1 \choose k+1}{2k+3 \choose k+2}{r-s \choose k-s+1}}{{2k+2 \choose k+1}{k+r+2 \choose
k+2}} \, a(R).
\end{eqnarray*}
Fix $R \supseteq S$, with $r > k+1$, and consider the coefficients of $a(R)$ in the previous equality. By equating them, we obtain an identity
which can be easily checked.\qed
\end{pf*}

\begin{exmp}\label{ex:1}
Let $\hat f : [0,1]^4 \to \R$ be given by
\begin{eqnarray*}
\hat f(\bx) &=& \frac 3{10} \, \bigl(x_1+x_2+x_3+\min\{x_1,x_2\}+\min\{x_1,x_3\}+\min\{x_2,x_3\}\bigr)\\
&& \mbox{}-\frac{21}{25} \, \min\{x_1,x_2,x_3\} + \frac 1{25} \, \min\{x_1,x_2,x_3,x_4\}.
\end{eqnarray*}
The best constant approximation is given by
$$
(A_0\hat f)(\bx) = \frac{137}{250} \, ,
$$
the best linear approximation by
$$
(A_1\hat f)(\bx) = \frac 1{100} + \frac{89}{250} \, (x_1+x_2+x_3) + \frac 1{125} \, x_4
$$
and the best min-quadratic approximation by
\begin{eqnarray*}
(A_2\hat f)(\bx) &=& -\frac{27}{700} + \frac{803}{1750} \, (x_1+x_2+x_3) -\frac{8}{875} \, x_4 \\
&& \null-\frac{19}{175}\, \bigl(\min\{x_1,x_2\}+\min\{x_1,x_3\}+\min\{x_2,x_3\}\bigr) \\
&& \null+\frac 2{175}\, \bigl(\min\{x_1,x_4\}+\min\{x_2,x_4\}+\min\{x_3,x_4\}\bigr).
\end{eqnarray*}
\end{exmp}

Before closing this section, we show that every approximation preserves the symmetry properties of the functions. For instance, in the previous
example, we observe that the function $\hat f$ and all its approximations are symmetric in the variables $x_1$, $x_2$, and $x_3$.

Let us establish this result in a more general setting. Let $X$ be a nonempty closed convex set in a finite-dimensional inner product space $V$.
For any $u\in V$, the distance between $u$ and $X$ is achieved at a unique point in $X$. We denote this point by $A_X(u)$ and call it the
projection of $u$ onto $X$ (see for instance \cite[Chap.~3, \S 3.1]{HirLem93}).

\begin{lem}\label{lemma:iac}
Let $I:V\to V$ be an isometry such that $I(X)=X$. Then $I$ and $A_X$ commute, i.e., $I\circ A_X=A_X\circ I$.
\end{lem}

\begin{pf*}{Proof.}
We clearly have $I[A_X(u)]\in X$ for any $u\in V$. Furthermore, for any $v\in X$, we have
$$
\|I(u)-I[A_X(u)]\| = \|u-A_X(u)\| \leqslant \|u-v\| = \|I(u)-I(v)\|,
$$
which shows that $I[A_X(u)]$ is the projection of $I(u)$ onto $X$.\qed
\end{pf*}

Let us now apply this result to our least squares approximation problem.

\begin{defn}
For any $\sigma$ in $\mathfrak{S}_n$, we define the linear operator $P_{\sigma}$ of $V_n$ as
\[
P_{\sigma}\hat f(x_{1},\ldots,x_{n}):=\hat f(x_{\sigma(1)},\ldots,x_{\sigma(n)}).
\]
We say that $\sigma\in\mathfrak{S}_n$ is a {\em symmetry}\/ of $\hat f$ if $P_{\sigma}(\hat f)=\hat f$.
\end{defn}

\begin{prop}
For any $\sigma\in\mathfrak{S}_n$, the operator $P_{\sigma}$ is an isometry of $V_n$ such that $P_{\sigma}(V_k)=V_k$. In particular, the
operators $P_{\sigma}$ and $A_k$ commute and $A_k$ preserves the symmetries of its arguments.
\end{prop}

\begin{pf*}{Proof.}
It is clear that $P_{\sigma}(V_k)=V_k$ and that $P_{\sigma}$ is an isometry, i.e., it satisfies
$$
\langle P_{\sigma} \hat f, P_{\sigma} \hat g\rangle = \langle \hat f, \hat g\rangle\qquad (\hat f, \hat g\in\,V_n).
$$
Therefore, by Lemma~\ref{lemma:iac}, $P_{\sigma}$ and $A_k$ commute. In particular, if $\sigma\in\mathfrak{S}_n$ is a symmetry of $\hat f$, then
$$P_{\sigma}[A_k(\hat f)]=A_k[P_{\sigma}(\hat f)]=A_k(\hat f)$$ and hence $\sigma$ is also a symmetry of $A_k(\hat f)$.\qed
\end{pf*}

\section{A new interaction index}

In cooperative game theory the concept of {\em power index}\/ (or {\em value}) was introduced in the pioneering work of Shapley \cite{Sha53}.
Roughly speaking, a power index on $N$ is a function $\phi:\mathcal{F}_N\times N\to\R$ that assigns to every player $i \in N$ in a game $v \in
\mathcal{F}_N$ his/her prospect $\phi(v,i)$ from playing the game. The {\em Shapley power index}\/ of a player $i \in N$ in a game $v\in
\mathcal{F}_N$ is given by
$$
\phi_{Sh}(v,i) := \sum_{T \subseteq N \setminus \{i\}}  \frac 1n{n-1\choose t}^{-1}[v(T \cup \{i\}) - v(T)].
$$

Another frequently used power index is the {\em Banzhaf power index}\/ \cite{Ban65,DubSha79} which, for a player $i \in N$ in a game $v \in
\mathcal{F}_N$, is defined by
$$
\phi_B(v,i) := \sum_{T \subseteq N \setminus \{i\}}  \frac{1}{2^{n-1}}\, [v(T \cup \{i\}) - v(T)].
$$

The concept of {\em interaction index}, which is an extension of that of power index, was recently introduced axiomatically to measure the
interaction phenomena among players. An interaction index on $N$ is essentially a function $I:\mathcal{F}_N\times 2^N\to\R$ that assigns to
every coalition $S\subseteq N$ of players in a game $v \in \mathcal{F}_N$ its interaction degree. Various interaction indices have been
introduced thus far in the literature: the {\em Shapley interaction index}\/ \cite{Gra97}, the {\em Banzhaf interaction index}\/
\cite{GraRou99,Rou96}, and the {\em chaining interaction index}\/ \cite{MarRou99}, which all belong to the class of {\em cardinal-probabilistic
interaction indices}\/ (see Definition~\ref{de:ip} below) newly axiomatized in \cite{FujKojMar06}.


For instance, the {\em Banzhaf interaction index}\/ on $N$, which extends the concept of Banzhaf power index on $N$, is the mapping
$I_{\mathrm{B}}:\mathcal{F}_N\times 2^N\to\R$ defined by
$$
I_{\mathrm{B}}(v,S):=\sum_{T\subseteq N\setminus S} \frac 1{2^{n-s}}\,\Delta_S v(T),
$$
where $\Delta_S v(T)$ is the $S$-derivative of $v$ at $T$ defined for any disjoint subsets $S,T\subseteq N$ by
\begin{equation}\label{eq:DiscDer}
\Delta_S v(T):=\sum_{R\subseteq S} (-1)^{s-r}\, v(R\cup T).
\end{equation}
This index can be easily expressed in terms of the M\"obius transform (\ref{eq:mo}) of $v$ as
$$
I_{\mathrm{B}}(v,S)=\sum_{T\supseteq S} \frac 1{2^{t-s}}\, a(T).
$$

It is noteworthy that, besides the axiomatic approach presented in \cite{GraRou99,Rou96}, the Banzhaf interaction index can also be defined from
the Hammer-Holzman approximation problem. Indeed, as pointed out in \cite{HamHol92} and \cite[\S 7]{GraMarRou00}, by considering the leading
coefficients in (\ref{eq:gr/ma/ro}) of the best $s$th approximation, for all $s=0,\ldots,n$, we immediately observe that
$$
I_{\mathrm{B}}(v,S) = a_s(S).
$$

In this section we use the same approach to define a new interaction index from our approximation problem of Lov\'asz extensions. In this sense
this new index can be seen as an analog of the Banzhaf interaction index.

Before going on, we recall the concept of cardinal-probabilistic interaction index \cite{FujKojMar06}.

\begin{defn}\label{de:ip}
A {\em cardinal-probabilistic interaction index}\/ on $N$ is a mapping $I_p:\mathcal{F}_N\times 2^N\to\R$ such that, for any $S\subseteq N$,
there is a family of nonnegative real numbers $\{p_t^s(n)\}_{t=0,\ldots,n-s}$ satisfying $\sum_{t=0}^{n-s}p_t^s(n)=1$, such that
\begin{equation}\label{eq:ip}
I_p(v,S)=\sum_{T\subseteq N\setminus S} p_t^s(n)\,\Delta_S v(T).
\end{equation}
\end{defn}

It has been proved \cite[\S4]{FujKojMar06} that a mapping $I_p:\mathcal{F}_N\times 2^N\to\R$ of the form (\ref{eq:ip}) is a
cardinal-probabilistic interaction index on $N$ if and only if, for any integer $s\in\{0,\ldots,n\}$, there exists a uniquely determined
cumulative distribution function $F_s$ on $[0,1]$ such that
\begin{equation}\label{eq:ptsn}
p_t^s(n)=\int_0^1 x^t(1-x)^{n-s-t}\, dF_s(x).
\end{equation}
Moreover, defining $q_t^s:=p_{t-s}^s(t)=\int_0^1 x^{t-s}\, dF_s(x)$ for all integers $s,t$ such that $0\leqslant s\leqslant t\leqslant n$, we
have
\begin{equation}\label{eq:ipa}
I_p(v,S)=\sum_{T\supseteq S} q_t^s\, a(T).
\end{equation}

The following lemma provides conditions on arbitrary coefficients $q_t^s$ so that a mapping of the form (\ref{eq:ipa}) can be a
cardinal-probabilistic interaction index.

\begin{lem}\label{lemma:qts}
Consider an infinite sequence $\{q_t^s\}_{t\geqslant s\geqslant 0}$. Then the following three conditions are equivalent:
\begin{itemize}
\item[(i)] we have $\sum_{i=0}^m(-1)^i{m\choose i}q_{t+i}^s\geqslant 0$ for all $m\geqslant 0$ and all $t\geqslant s\geqslant 0$,

\item[(ii)] there is a unique cumulative distribution function $F_s$ on $[0,1]$ such that, for any integer $t\geqslant s\geqslant 0$,
$$
q_t^s=\int_0^1 x^{t-s} \, dF_s(x),
$$

\item[(iii)] for any $N=\{1,\ldots,n\}$, the function $I_p:\mathcal{F}_N\times 2^N\to\R$, defined in (\ref{eq:ipa}) is a cardinal-probabilistic
interaction index.
\end{itemize}
Moreover, when these conditions are satisfied, then the mapping $I_p:\mathcal{F}_N\times 2^N\to\R$, defined in (\ref{eq:ipa}), is of the form
(\ref{eq:ip}) and the coefficients $q_t^s$ and $p_t^s(n)$ are linked through $q_t^s=p_{t-s}^s(t)$ and
$$
p_t^s(n)=\sum_{i=s+t}^n (-1)^{i-s-t}{n-s-t\choose i-s-t}q_i^s.
$$
\end{lem}

\begin{pf*}{Proof.}
$(i)\Leftrightarrow (ii)$ Follows immediately from the Hausdorff's moment problem (see for instance Akhiezer \cite[Theorem 2.6.4]{Akh65}).

$(ii)\Rightarrow (iii)$ Rewriting (\ref{eq:ipa}) in terms of $v$ gives
$$
I_p(v,S) = \sum_{J\supseteq S}q_j^s\, \sum_{L\subseteq J}(-1)^{j-l} v(L).
$$
Partitioning $L\subseteq J$ into $R\subseteq S$ and $T\subseteq J\setminus S$ and then using (\ref{eq:DiscDer}), we obtain
\begin{eqnarray*}
I_p(v,S) &=& \sum_{T\subseteq N\setminus S}\bigg[\sum_{J\supseteq S\cup T}
(-1)^{j-s-t}q_j^s\bigg]\sum_{R\subseteq S}(-1)^{s-r}\, v(R\cup T)\\
&=& \sum_{T\subseteq N\setminus S}p_t^s(n)\Delta_S v(T),
\end{eqnarray*}
where the coefficients $p_t^s(n)$ are given by
$$
p_t^s(n)=\sum_{j=s+t}^n (-1)^{j-s-t}{n-s-t\choose j-s-t}q_j^s = \int_0^1 x^t(1-x)^{n-s-t}\, dF_s(x),
$$
which shows that $I_p$ is a cardinal-probabilistic interaction index.

$(iii)\Rightarrow (ii)$ According to formula (\ref{eq:ipa}), we have
$$
I_p(v,S)=\sum_{T\supseteq S} q_t^s\, a(T)=\sum_{T\supseteq S} \Big[\int_0^1 x^{t-s} \, dF_s(x)\Big]\, a(T),
$$
which completes the proof by uniqueness of this decomposition.\qed
\end{pf*}


We are now ready to introduce our new interaction index, namely the mapping
$$
I_M:\mathcal{F}_N\times 2^N\to\R\, :\, (v,S)\mapsto I_M(v,S)=a_s(S),
$$
where $a_s(S)$ are the leading coefficients in (\ref{eq:solution/system}) of the best $s$th approximation of $\hat f$, for all $s=0,\ldots,n$.

This immediately leads to the following equivalent definition, which makes use of the classical beta function
$$
B(a,b):=\int_0^1u^{a-1}(1-u)^{b-1}\, du \qquad (a,b>0).
$$

\begin{defn}
For any $N=\{1,\ldots,n\}$, we define the mapping $I_M:\mathcal{F}_N\times 2^N\to\R$ as
$$
I_M(v,S)=\sum_{T\supseteq S} q_t^s\, a(T),
$$
where $q_t^s:=\frac{{2s+1\choose s+1}}{{s+t+1\choose s+1}}=\frac{B(t+1,s+1)}{B(s+1,s+1)}$ for all $t\geqslant s\geqslant 0$.
\end{defn}

By using Lemma~\ref{lemma:qts}, we can easily show that $I_M$ is a cardinal-probabilistic interaction index.

\begin{prop}\label{prop:ims}
The function $I_M:\mathcal{F}_N\times 2^N\to\R$ is a cardinal-probabilistic interaction index, given by
$$
I_M(v,S)=\sum_{T\subseteq N\setminus S} p_t^s(n)\, \Delta_S\, v(T),
$$
where $p_t^s(n):=\frac{B(n-t+1,s+t+1)}{B(s+1,s+1)}$.
\end{prop}

\begin{pf*}{Proof.}
We immediately observe that Lemma~\ref{lemma:qts} applies with the beta distribution
$$
F_s(x)=\frac{\int_0^x u^s(1-u)^s\, du}{B(s+1,s+1)}.
$$
Indeed, we have
$$
\int_0^1 x^{t-s} \, dF_s(x)=\frac{1}{B(s+1,s+1)}\int_0^1 x^{t}(1-x)^s \, dx=q_t^s.
$$
The expression of the coefficients $p_t^s(n)$ then follows immediately from (\ref{eq:ptsn}).\qed
\end{pf*}

The following proposition shows that the mapping $v\mapsto (I_M(v,S))_{S\subseteq N}$ is invertible and yields the transformation formula from
$(I_M(v,S))_{S\subseteq N}$ to $a=(a(S))_{S\subseteq N}$.

\begin{prop}
For any $S\subseteq N$, we have
$$
a(S)=\sum_{T\supseteq S} h_t^s\, I_M(v,T).
$$
with
$$
h_t^s =(-1)^{t-s}\,\frac{{s+t\choose t}}{{2t\choose t}}.
$$
\end{prop}

\begin{pf*}{Proof.}
We only need to show that
$$
\sum_{T\supseteq S} q_t^s\,\sum_{R\supseteq T} h_r^t\ I_M(v,R) = I_M(v,S).
$$
We have
$$
\sum_{T\supseteq S} q_t^s\,\sum_{R\supseteq T} h_r^t\ I_M(v,R) = \sum_{R\supseteq S} I_M(v,R)\, \sum_{T:\, S\subseteq T\subseteq R} q_t^s\,
h_r^t,
$$
where the inner sum equals 1 if $S=R$ and 0 if $S\varsubsetneq R$. Indeed, in the latter case, we have
\begin{eqnarray*}
\sum_{T:\, S\subseteq T\subseteq R} q_t^s\, h_r^t
&=& \frac{{2s+1\choose s+1}}{{2r\choose r}}\,\frac{(s+1)!}{r!}\,\sum_{t=s}^r (-1)^{r-t} {r-s\choose t-s}\frac{(r+t)!}{(s+t+1)!}\\
&=& \frac{{2s+1\choose s+1}}{{2r\choose r}}\,\frac{(s+1)!}{r!}\,\sum_{t=s}^r (-1)^{r-t} {r-s\choose t-s}\Big[\frac{d^{r-s-1}}{dx^{r-s-1}}\,
x^{r+t}\Big]_{x=1}\\
&=& \frac{{2s+1\choose s+1}}{{2r\choose r}}\,\frac{(s+1)!}{r!}\,\Big[\frac{d^{r-s-1}}{dx^{r-s-1}}\, x^{r+s}(x-1)^{r-s}\Big]_{x=1}\\
&=& 0.\qed
\end{eqnarray*}
\end{pf*}

In order to conclude, we focus on the power index associated to $I_M$. The following corollary shows that, incidentally, this power index
identifies with the index $W_i(\hat f)$ previously introduced axiomatically by Grabisch and Labreuche \cite[Theorem 2]{GraLab01} in the context
of multicriteria decision making.

\begin{cor}
The restriction of $I_M(v,\cdot)$ to singletons is given by:
$$
I_M(v,\{i\})= \sum_{T\subseteq N\setminus \{i\}} \frac{6(n-t)!(t+1)!}{(n+2)!}\,[v(T\cup\{i\})-v(T)].
$$
\end{cor}

\section*{Appendix A: Proof of Lemma \ref{lemma:integral}}

Observe first that we can assume that $S$ and $T$ are such that $|S \cup T| = n$.
Moreover, suppose that $S$ and $T$ are nonempty. Remark that $[0,1]^n$  is equal almost everywhere to the disjoint union of the sets
\[
S_{\sigma}^{\circ} := \{\bx \in [0,1]^n \mid x_{\sigma(1)} < \cdots < x_{\sigma(n)}\} \qquad (\sigma \in \mathfrak{S}_n).
\]
Then we have
\begin{equation}\label{sum}
\int_{[0,1]^n} {\rm min}_S(\bx) \, {\rm min}_T(\bx)\, d\bx=\sum_{\sigma \in \mathfrak{S}_n}I_\sigma,
\end{equation}
where
\[
I_\sigma := \int_{S_{\sigma}^{\circ}} {\rm min}_S(\bx) \, {\rm min}_T(\bx)\, d\bx\qquad (\sigma\in \mathfrak{S}_n),
\]

Now, define for every $R \subseteq N$ and $p\in\{1,\ldots,|R|\}$ the set
\[\mathfrak{S}_n^{(p)}(R) := \{\sigma \in\mathfrak{S}_n \mid
\sigma(1),\ldots,\sigma(p) \in R,\,\sigma(p+1)\not\in R\}.\] It is then clear that there is a decomposition of $\mathfrak{S}_n$ into disjoint
subsets
\begin{equation}\label{eq:disjunion}
\mathfrak{S}_n = \Bigl[\mathfrak{S}_n(S\cap T)\Bigr] \cup \Bigl[\, \bigcup_{p=1}^{n-t} \mathfrak{S}_n^{(p)}(S \setminus T)\Bigr] \cup \Bigl[\,
\bigcup_{p=1}^{n-s} \mathfrak{S}_n^{(p)}(T \setminus S)\Bigr],
\end{equation}
where
\[\mathfrak{S}_n(S\cap T) := \{\sigma \in\mathfrak{S}_n \mid
\sigma(1)\in S\cap T\}.\] Moreover we have
\[\left\{\begin{array}{lll}
|\mathfrak{S}_n(S \cap T)| = |S \cap T| \, (n-1)!\\
 |\mathfrak{S}_n^{(p)}(S \setminus T)| = {n-t \choose p} \, p! \, t \,
 (n-p-1)!\, .
\end{array}\right.\]
Indeed, every element of $\mathfrak{S}_n(S \cap T)$ corresponds to a
 choice of $\sigma(1)$ in $S \cap T$, and a permutation of the remaining
 elements in $N$. In the same fashion, every element of
 $\mathfrak{S}_n^{(p)}(S \setminus T)$ is determined by the choice of $p$
 elements in $S\setminus T$, a permutation of these elements, the choice of
 a single element in $T$ and a permutation of the remaining elements.
Moreover, for any $\sigma \in \mathfrak{S}_n(S \cap T)$, we have
$$
I_\sigma = \int_0^1 \, \int_0^{x_{\sigma(n)}} \cdots \int_0^{x_{\sigma(2)}} x_{\sigma(1)}^2
 \, dx_{\sigma(1)} \cdots dx_{\sigma(n)} =
\frac 2{(n+2)!}\, ,
$$
and
\begin{equation}\label{eq:sum1}
\sum_{\sigma \in \mathfrak{S}_n(S \cap T)} I_\sigma = \frac{2 \, |S \cap T| \, (n-1)!} {(n+2)!} = \frac{2 \, (s+t-n)}{n(n+1)(n+2)}\, .
\end{equation}
In the same way, for any $\sigma \in \mathfrak{S}_n^{(p)}(S \setminus T)$, we have
\begin{eqnarray*}
I_\sigma &=& \int_0^1 \, \int_0^{x_{\sigma(n)}} \cdots \int_0^{x_{\sigma(p+2)}} x_{\sigma (p+1)} \int_0^{x_{\sigma(p+1)}} \cdots
\int_0^{x_{\sigma(2)}} x_{\sigma(1)} \, dx_{\sigma(1)} \cdots dx_{\sigma(n)} \\
&=& \frac {p+2}{(n+2)!}\, ,
\end{eqnarray*}
and
\begin{eqnarray*}
\sum_{p=1}^{n-t} \, \sum_{\sigma \in \mathfrak{S}_n^{(p)}(S \setminus T)} I_\sigma & = & \frac{t}{(n+2)!} \sum_{p=1}^{n-t} {n-t \choose p} \, p!
\, (n-p-1)! \,
(p+2)\\
& = & \frac{t! \, (n-t)!}{(n+2)!} \, \sum_{p=1}^{n-t} {n-p-1 \choose t-1} (p+2)
\\
& = & \frac{t! \, (n-t)!}{(n+2)!} \, \biggl[(n+2) \, \sum_{p=1}^{n-t}{n-p-1
\choose t-1} - t \, \sum_{p=1}^{n-t}{n-p \choose t}\biggr]\\
& = & \frac{t! \, (n-t)!}{(n+2)!} \, \biggl[(n+2){n-1 \choose t} - t {n \choose t+1}\biggr]\, ,
\end{eqnarray*}
that is,
\begin{equation}\label{eq:sum2}
\sum_{p=1}^{n-t} \, \sum_{\sigma \in \mathfrak{S}_n^{(p)}(S \setminus T)} I_\sigma = \frac{(n-t)(n+2t+2)}{n(n+1)(n+2)(t+1)}\, .
\end{equation}
Similarly, we can write
\begin{equation}\label{eq:sum3}
\sum_{p=1}^{n-s} \, \sum_{\sigma \in \mathfrak{S}_n^{(p)}(T \setminus S)} I_\sigma = \frac{(n-s)(n+2s+2)}{n(n+1)(n+2)(s+1)}
\end{equation}
and by Eqs.\ (\ref{sum})--(\ref{eq:sum3}), we have
\begin{eqnarray*}
\int_{[0,1]^n} {\rm min}_S(\bx)\,{\rm min}_T(\bx) \, d\bx & = & \sum_{\sigma \in \mathfrak{S}_n(S \cap T)} I_\sigma + \sum_{p=1}^{n-t} \,
\sum_{\sigma \in \mathfrak{S}_n^{(p)}(S \setminus T)} I_\sigma + \sum_{p=1}^{n-s} \, \sum_{\sigma \in
\mathfrak{S}_n^{(p)}(T \setminus S)} I_\sigma \\
& = & \frac{s+t+2}{(n+2)(s+1)(t+1)}
\end{eqnarray*}
as desired. We can easily see that the result still holds if $S = \varnothing$
 or $T = \varnothing$.\qed

\section*{Appendix B: Proof of Lemma \ref{lemma:enfer/1}}

We only have to prove that the explicit expression we give for $A_k({\rm min}_{S})$ fulfills the conditions of Proposition \ref{exist}, namely,
$$
\big\langle {\rm min}_{S}- A_k({\rm min}_{S}),\,{\rm min}_{T}\big\rangle=0 \qquad (T\subseteq N,\, t \leqslant k),
$$
or, equivalently,
\begin{equation}\label{formule}
\sum_{J \subseteq S} I_{J,T} \, (-1)^{k+j} \, {k+j+1 \choose k+1} = 0 \qquad (T\subseteq N,\, t \leqslant k),
\end{equation}
where $I_{J,T}:=\langle {\rm min}_J, {\rm min}_T\rangle$ is given explicitly in (\ref{eq:int/min/min}).

Partitioning $J \subseteq S$ into $P \subseteq S \setminus T$ and $Q \subseteq R$ with $R:=S \cap T$ (hence $r\leqslant t\leqslant k$), the
left-hand side of (\ref{formule}) becomes
\begin{eqnarray*}
\lefteqn{\sum_{P \subseteq S \setminus T} \: \sum_{Q \subseteq R} I_{P \cup Q,T} \, (-1)^{k+p+q} {k+p+q+1 \choose k+1}}\\
&=& \sum_{p=0}^{k-r+1} {k-r+1 \choose p} \sum_{q=0}^r {r \choose q} \, \frac{p+q+t+2}{(p+t+2)(p+q+1)(t+1)} \, (-1)^{k+p+q} {k+p+q+1 \choose
k+1}\\
&=& - \sum_{p=0}^{k-r+1} {k-r+1 \choose p} (-1)^{k-r+1-p} \, \frac 1{p+t+2} \, \sum_{q=0}^r {r \choose q} (-1)^{r-q} \, g(p+q)
\end{eqnarray*}
where $g : \N \to \R$ is defined by
\begin{eqnarray*}
g(z) & = & \frac{z+t+2}{(z+1)(t+1)} {z+k+1 \choose k+1}\\
& = & \frac 1{k+1} {z+k+1 \choose k} + \frac 1{t+1} {z+k+1 \choose k+1}.
\end{eqnarray*}

Consider the difference operator
$$
\Delta_n \, f(n) := f(n+1)-f(n)
$$
for functions on $\N$. It is well-known that we have
\begin{equation}\label{eq:diff}
\Delta_n^k \, f(n) = \sum_{j=0}^k {k \choose j} (-1)^{k-j} \, f(n+j) \qquad (k \in \N).
\end{equation}
Applying this to $g$, we obtain
$$
\sum_{q=0}^r {r \choose q} (-1)^{r-q} \, g(p+q) = \Delta_p^r \, g(p) = \frac 1{k+1} {p+k+1 \choose k-r} + \frac 1{t+1} {p+k+1 \choose k-r+1}
$$
and the left-hand side of (\ref{formule}) becomes
$$
- \sum_{p=0}^{k-r+1} {k-r+1 \choose p} (-1)^{k-r+1-p} \, \frac 1{p+t+2} \, \Bigl(\frac 1{k+1} {p+k+1 \choose k-r} + \frac 1{t+1} {p+k+1 \choose
k-r+1}\Bigr)\, .
$$
Applying (\ref{eq:diff}) again, we see that this latter expression can be written as
$$
- \Bigg[ \Delta_z^{k-r+1} \, \frac 1{z+t+2} \, \Bigl(\frac 1{k+1} {z+k+1 \choose k-r} + \frac 1{t+1} {z+k+1 \choose k-r+1}\Bigr)\Bigg]_{z = 0}.
$$
We now show that the expression in brackets is identically zero. We do this by considering two cases:
\begin{itemize}
\item If $t = k$ then
\begin{eqnarray*}
\lefteqn{\frac 1{z+t+2} \, \Bigl(\frac 1{k+1} {z+k+1 \choose k-r} + \frac 1{t+1} {z+k+1 \choose k-r+1}\Bigr)}\\
& = & \frac 1{(k+1)(z+k+2)} \, {z+k+2 \choose k-r+1} \enskip = \enskip  \frac 1{(k+1)(k-r+1)} \, {z+k+1 \choose k-r}
\end{eqnarray*}
is a polynomial $P_{k-r}(z)$ of degree $k-r$, and $\Delta_z^{k-r+1} \, P_{k-r}(z) \equiv 0$.

\item If $r \leqslant t \leqslant k-1$ then
\begin{eqnarray*}
\lefteqn{\frac 1{z+t+2} \, \Bigl(\frac 1{k+1} {z+k+1 \choose k-r} + \frac 1{t+1} {z+k+1 \choose k-r+1}\Bigr)}\\
& = & \frac 1{z+t+2} \, \Bigl(\frac 1{(k+1)(k-r)!} \prod_{i=r+2}^{k+1}(z+i) + \frac 1{(t+1)(k-r+1)!} \prod_{i=r+1}^{k+1}(z+i)\Bigr)
\end{eqnarray*}
is a polynomial $Q_{k-r}(z)$ of degree $k-r$, and $\Delta_z^{k-r+1} \, Q_{k-r}(z) \equiv 0$.\qed
\end{itemize}


\end{document}